\RequirePackage{fix-cm}
\documentclass{article}
\usepackage{amsmath}
\usepackage{amssymb}
\usepackage{indentfirst}
\usepackage{graphicx}
\usepackage{epstopdf}
\usepackage{caption}

\date{}

\numberwithin{equation}{section}

\usepackage{amssymb,amsthm,amsmath}
\usepackage[numbers,sort&compress]{natbib}

\begin{document}

\newtheorem{theorem}{Theorem}[section]
\newtheorem{lemma}{Lemma}[section]
\newtheorem{definition}{Definition}[section]
\newtheorem{remark}{Remark}[section]
\newtheorem{example}{Example}[section]
\newtheorem{algorithm}{Algorithm}[section]
\setcounter{section}{1}
\baselineskip 16pt \setcounter{section}{0}

\title{\bf{\  A note on an Counter Example of Dacorogna}
\thanks{Supported by the National Natural Science Foundation of China (No.11671278)
}}
 \author{Yan Tang$^1$ and Shiqing Zhang$^2$\thanks{Corresponding author Shiqing Zhang, E-mail:
 zhangshiqing@msn.com}\\
{\small 1.College of  Mathematics and Statistics,Chongqing Technology and Business University,Chongqing,400067, China}\\
{\small 2. Department of  Mathematics, Sichuan
University,Chengdu,610064, China} }
\maketitle{}

\begin{quote}{\bf Abstract:}  In this short note,we correct a well-known counter example of the famous book of Dacorogna[2].

\end{quote}
\begin{quote}
\end{quote}
\begin{quote}
{\bf 2000 MR Subject Classification}:\: 49J10;49J45;74B20
\end{quote}
\section{Introduction}

It is well-known that all finite dimensional normed spaces are reflexive Banach spaces,any bounded closed subset in a finite dimensional normed space is compact,but for infinite dimensional normed spaces,the situations are very complicated, we know that $L^p(\Omega)$ and $W^{1,p}(\Omega) (1<p<+\infty)$ are reflexive Banach spaces,but $C(\Omega)$ is not.\\

 In 1940 and 1947,Shimulyan and Eberlein studied the necessary and sufficient condition that a normed vecor space is a reflexive Banach space, they obtain the following important Theorem.

{\bf Theorem 1.1} (Eberlein-Shimulyan[3][5]) The  normed vector space $X$ is reflexive if and only if any bounded sequence in $X$ has a convergence subsequence.

 A well-known result for the existence of the minimizer for a functional is the following Theorem.

 {\bf Theorem 1.2}([1][2]) Let $X$ be a reflexive Banach space,$M\neq\emptyset$ be weakly closed subset of $X$, $f:M\rightarrow \mathbb{R}\cup\{+\infty\}$, $f\neq +\infty$, $\inf_{x\in M}f(x)>-\infty$. Furthermore, we assume \\ \indent
 $\mathbf{(i).}$ $f$ is coercive, that is,$f(x)\rightarrow +\infty$ as $\|x\|\rightarrow+\infty$;\\ \indent
 $\mathbf{(ii).}$ $f$ is weakly lower semi-continuous.\\
 Then $\exists \bar{x}$ such that $f(\bar{x})=\inf_{x\in M}f(x)$.

 The reflexive property of the Banach space $X$ in Theorem 1.2 is a key condition, in the pages of 48 and 49 of Dacorogna's book, he pointed out '\textit{ The hypothesis on the reflexivity of $X$ cannot be dropped in general}' and he gave an counter example (Professor Docorogna also pointed out the example can also refer to Rudin[4]), he let
  \begin{eqnarray*}
 X=\{u\in C[0,1]:\int_0^{1/2}u(t)dt-\int_{1/2}^1u(t)dt=1\},
\end{eqnarray*}
 and let
  \begin{eqnarray*}
 \inf\{I(u)=\|u\|_{L^\infty}:u\in X\}.
\end{eqnarray*}
He observed
\begin{eqnarray*}
I(u)=\|u\|_{L^\infty}\geq\int_0^1|u(t)|dt\geq\int_0^{1/2}u(t)dt-\int_{1/2}^1u(t)dt=1,
\end{eqnarray*}
and he got
 \begin{eqnarray*}
 \inf(P)=\inf\{\|u\|_{L^\infty}\}=1.
\end{eqnarray*}

By using a minimizing sequence
\begin{equation}
u_N(x)=\left\{
\begin{array}{ll}
1,&x\in [0,\frac{1}{2}-\frac{1}{N}];\\
-Nx+\frac{N}{2},&x \in[\frac{1}{2}-\frac{1}{N}, \frac{1}{2}+\frac{1}{N} ];\\
-1,&x\in[ \frac{1}{2}+\frac{1}{N},1],
 \end{array}\right.
\end{equation}
 and it is clear that no continuous function can satisfy
 \begin{eqnarray*}
\|u\|_{L^\infty}=\int_0^1|u(t)|dt\geq\int_0^{1/2}u(t)dt-\int_{1/2}^1u(t)dt=1.
\end{eqnarray*}
 So, he said the hypothesis on the reflexivity of $X$ in Theorem1.2 cannot be dropped in general.

\section{Corrections}
 In the above example, we found $u_N(x)\notin X$.

  Since
  \begin{eqnarray*}
\int_0^{1/2}u(t)dt-\int_{1/2}^1u(t)dt=1-\frac{1}{n}\neq 1.
\end{eqnarray*}
 On the other hand, it's easy to see  $\|u_N\|_{L^\infty}=1$,so if $ \inf(P)=\inf\{\|u\|_{L^\infty}\}=1$,then $I(u)= \|u\|_{L^\infty}$ attains its infimum at $u_N(x)$, which is also a contradiction.

 In the following, we try to revise the function $u_N(x)$. We define

 \begin{equation}
u_n(x)=\left\{
\begin{array}{ll}
1+\dfrac{8n}{(n-2)^2}x,&x\in [0,\dfrac{1}{4}-\dfrac{1}{2n}];\\
1+\dfrac{4}{n-2}-\dfrac{8n}{(n-2)^2}x,&x \in[\dfrac{1}{4}-\frac{1}{2n}, \dfrac{1}{2}-\frac{1}{n} ];\\
-nx+\dfrac{n}{2}, &x\in[\dfrac{1}{2}-\dfrac{1}{n},\dfrac{1}{2}+\dfrac{1}{n} ];\\
-1+\dfrac{4(n+2)}{(n-2)^2}-\dfrac{8n}{(n-2)^2}x, &x\in [\dfrac{1}{2}+\dfrac{1}{n},\dfrac{3}{4}+\frac{1}{2n} ];\\
-1-\dfrac{8n}{(n-2)^2}+\dfrac{8n}{(n-2)^2}x,&
x\in[\dfrac{3}{4}+\dfrac{1}{2n},1 ]
 \end{array}\right.
\end{equation}

Here we will insert the figure of the function $u_n(x)$.
{\scriptsize
\begin{center}
\begin{figure}[hbt]
  \centering
  \includegraphics[width=0.3\textwidth]{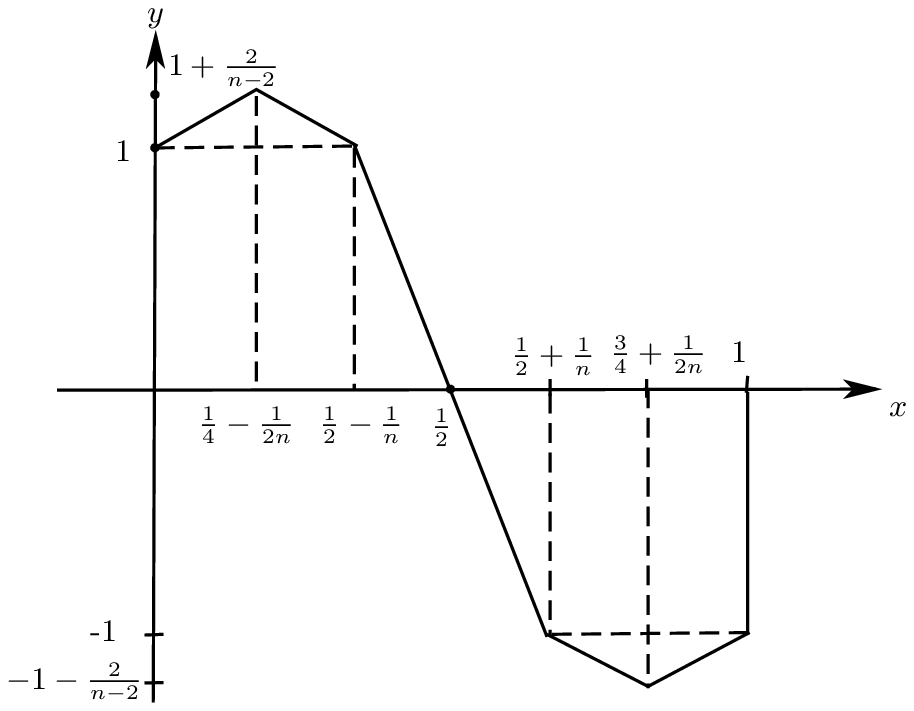}
  \vskip -0.2in
 \end{figure}
\end{center}}

Then it is not difficult to compute
$\int_0^{1/2}u(t)dt-\int_{1/2}^1u(t)dt=1$. Hence $u_{n}(x)\in X$,
furthermore,
  \begin{eqnarray*}
\|u\|_{L^\infty}=1+\frac{2}{n-2}\rightarrow 1.
\end{eqnarray*}
Consequently,
 \begin{eqnarray*}
 \inf(P)=\inf\{\|u\|_{L^\infty}\}=1.
\end{eqnarray*}
 It is easy to see that $X$ is  a convex strongly closed subset of $C[0,1]$, so it is also weakly closed subset of  $C[0,1]$. We know that the norm functional $\|u\|_{L^\infty} $ is weakly lower semi-continuous and coercive,but $C[0,1]$ is not a reflexive Banach space, so we can't apply Theorem1.2.

 In fact, it is easy to see that no continuous function $\tilde{u}$ satisfies
 \begin{eqnarray*}
\|\tilde{u}\|_{L^\infty}=\int_0^1|\tilde{u}(t)|dt=\int_0^{\frac{1}{2}}\tilde{u}(t)dt-\int_{\frac{1}{2}}^1\tilde{u}(t)dt=1.
\end{eqnarray*}

 \end{document}